\documentclass[12pt,oneside]{article}

\usepackage{amsfonts}
\usepackage{amscd}
\usepackage{amsthm}
\usepackage{amssymb}
\usepackage{amsmath}
\usepackage{epsfig}

\usepackage{graphicx}
\usepackage{xypic}
\input xypic
\input xy
\xyoption{all}

\newcommand{\mZ}{\mathbb Z}
\newcommand{\mR}{\mathbb R}
\newcommand{\mD}{\mathbb D}

\newcommand{\mG}{\Gamma}

\newcommand{\w}{\widehat}
\newcommand{\invlim}{\underleftarrow \lim}

\title{Involutions of negatively curved groups with wild boundary behavior.}

\author{F. T. Farrell \footnote{This research was supported in part
by the National Science Foundation.} \hskip 5pt \& J.-F. Lafont}

\theoremstyle{definition}

\theoremstyle{proposition}

\newtheorem{Prop}{Proposition}

\newtheorem{Thm}{Theorem}[section]
\newtheorem{Cor}{Corollary}[section]

\theoremstyle{plain}

\theoremstyle{remark}

\newtheorem*{Prf}{Proof}
\newtheorem{Rmk}{Remark}

\begin{document}

\maketitle

\begin{abstract}
We are interested in examples of compact, complete, locally CAT(-1) spaces $X$, 
and closed totally 
geodesic codimension two subspaces $Y$, with the property that $\partial
^\infty \tilde X=S^{n+2}_\infty$, and $\partial ^\infty \tilde Y=S^n_\infty$.  
We show that if the 
inclusion $S^n_\infty\hookrightarrow S^{n+2}_\infty$ induced by the inclusion 
$Y\subset X$ is a
non-trivial knot, then it is a totally wild knot (i.e. nowhere tame).  
We give examples where the inclusion $S^n_\infty\hookrightarrow S^{n+2}_\infty$ 
is indeed a 
non-trivial knot.
Furthermore, the examples $Y\subset X$ we construct can be realized as fixed
point sets of involutive isometries of $X$, so that the corresponding totally wildly 
knotted $S^n_\infty\hookrightarrow S^{n+2}_\infty$ are fixed point sets of geometric 
involutions
on $S^{n+2}_\infty=\partial ^\infty \tilde X$.  In the appendix, we also include a 
complete
criterion for knottedness of tame codimension two spheres in high dimensional 
($\geq 6$) spheres.
\end{abstract}

\section{Introduction.}

For a space $Z$ that arises as the boundary at infinity of a 
$\delta$-hyperbolic group $\mG$, it is 
of some interest to determine which self-homeomorphisms of $Z$
are {\it algebraic}, i.e. are induced by an 
automorphism of $\mG$.  Likewise, for a space $Z$ that arises as the 
boundary at infinity of the universal cover of a compact locally CAT(-1) 
space $X$, 
one can try to determine which self-homeomorphisms of $Z$ are {\it geometric},
i.e. are induced by isometric self-maps of $X$.  Note that geometric 
self-homeomorphisms induce algebraic self-homeomorphisms (by setting $\mG=\pi_1(X)$).
In a previous paper [14], the authors constructed geometric involutions of 
spheres whose fixed point set had 
infinitely generated (\v{C}ech) cohomology.  In this note, we analyze the
geometry of fixed point sets of geometric involutions in the special case
where $Z=S^{n+2}$ and the fixed point sets 
are embedded codimension two spheres.

Before stating our results, we remind the reader of a few basic facts
about embeddings $S^n\hookrightarrow S^{n+2}$.  For any such embedding,
Alexander Duality tells us that the space $S^{n+2}-S^n$ is a homology
$S^1$, hence there is a canonical infinite cyclic cover $\w{S^{n+2}-S^n}
\rightarrow S^{n+2}-S^n$ arising from the abelianization $\pi_1(S^{n+2}-S^n)
\twoheadrightarrow H_1(S^{n+2}-S^n)\cong \mZ$.  An embedding 
$S^n\hookrightarrow 
S^{n+2}$ is {\it unknotted} if the complement $S^{n+2}-S^n$ is homeomorphic
to $S^1\times \mR^{n+1}$, otherwise it is said to be {\it knotted}.  Note
that, for an unknotted embedding, the infinite cyclic cover $\w{S^{n+2}-S^n}$ is 
homeomorphic to $\mR\times \mR^{n+1}$, while the fundamental group of the
complement is $\pi_1(S^{n+2}-S^n)= \mZ$.  This gives sufficient 
conditions for an embedding $S^n\hookrightarrow S^{n+2}$ to be knotted, namely
if either:
\begin{itemize}
\item $\pi_1(S^{n+2}-S^n)\neq \mZ$,
\item or $H_s(\w{S^{n+2}-S^n})\neq 0$ for some $s\neq 0,n$,
\end{itemize}
then the embedding is knotted.  Conversely, when $n\geq 4$ and the embedding 
is smooth,  one can use work of Browder-Levine [7] to show that the above 
conditions are also necessary (see Theorem 4.1 in the Appendix for a proof).

An embedding $S^n\hookrightarrow S^{n+2}$ is said to be {\it tame} at a 
point $p\in S^n$ provided there exists a neighborhood $U$ of $p$ in 
$S^{n+2}$ with the property that $(U,S^n\cap U)$ is homeomorphic to 
$(\mR^{n+2}, \mR^n)$ (where the $\mR^n\hookrightarrow \mR^{n+2}$ is the 
embedding into the first factor of $\mR^{n+2}\cong \mR^n\times \mR^2$).
In other words, the embedding is ``collarable'' near the point $p$.  A
point $p\in S^n$ that is not tame is said to be {\it wild}.  
An embedding is tame provided it is tame at every point, and is 
{\it totally wild} if it is wild at every point.  Examples of
non-tame embeddings of $2$-spheres in $3$-spheres 
include the Alexander horned sphere, the Antoine wild sphere and the 
Fox-Artin ball (see Bing [6] for details).  A result which will be 
frequently used is that codimension two tame embeddings have tubular neighborhoods; 
this is due to Kirby-Siebenmann [17] when the ambient dimension is $\neq 4$,
and to Freedman-Quinn (Section 9.3 in [15]) in the remaining case.  We are now ready 
to state our main results. 

\begin{Thm}
Let $\Gamma$ be a $\delta$-hyperbolic group, $\Lambda\leq \Gamma$ a 
quasi-convex subgroup.  Assume that $\partial^\infty \Gamma =S^{n+2}_\infty$, $\partial
^\infty \Lambda=S^n_\infty$, where $n\neq 0,3$, and that the embedding $S^n_\infty 
\hookrightarrow S^{n+2}_\infty$ induced by the inclusion $\Lambda\leq \Gamma$ is 
knotted.  Then the embedding is a totally wild knot.
\end{Thm}

This allows us to answer a question we asked in our previous paper [14]:

\begin{Cor}[Generalized Smith conjecture for algebraic involutions]
Let $\tau:S^{n+2}\rightarrow S^{n+2}$ be an algebraic self-homeomorphisms 
($n\neq 0,3$) of finite order, and assume that
the fixed point set is an embedded $S^n\hookrightarrow S^{n+2}$.  If the embedding
is tame, then $S^n\hookrightarrow S^{n+2}$ is unknotted.
\end{Cor}

In the previous corollary, by an algebraic self-homeomorphism of finite order, 
we mean an algebraic self-homeomorphism induced by an automorphism of finite 
order.  Theorem 1.1 also yields the immediate corollary:

\begin{Cor}
Let $X$ be a compact, complete, locally CAT(-1) space, $Y\subset X$ a totally
geodesic subspace, $\tilde Y$ and $\tilde X$ their respective universal covers.
Assume that $\partial ^\infty \tilde Y=S^n_\infty$, $\partial ^\infty \tilde X=
S^{n+2}_\infty$, where $n\neq 0,3$,
and that the embedding $S^n_\infty\hookrightarrow S^{n+2}_\infty$ induced
by a lift $\tilde Y\subset \tilde X$ is knotted.  Then the embedding is a 
totally wild knot.
\end{Cor}

Even in the special case of fundamental groups of negatively curved Riemannian 
manifolds, this gives the interesting:

\begin{Cor}
Let $M^n$ ($n\geq 4$, $n\neq 6$) be a closed, negatively curved Riemannian manifold, 
and let $\alpha\in
Aut(\pi_1(M^n))$ have finite order.  Denote by $\alpha_\infty:S^{n-1}_\infty
\rightarrow S^{n-1}_\infty$ the induced involution on the sphere at infinity
of the universal cover $\tilde M^n$.  Then the fixed point set of $\alpha_\infty$
cannot be a tame codimension two knot.
\end{Cor}

The next reasonable question is whether one can find non-trivial knots.  This
is addressed in the following:

\begin{Thm}
Let $S^n\hookrightarrow S^{n+2}$ be 
a tamely knotted codimension two PL-embedding, with $n\geq 5$.  
Let $X$ be the strict hyperbolization
of the suspension of $S^{n+2}$, and $Y\subset X$ the strict hyperbolization of
the suspension of $S^n$.  Then the embedding $\partial ^\infty \tilde Y
\subset \partial ^\infty \tilde X$ is a 
knotted $S^n_\infty$ in $S^{n+2}_\infty$.  By the previous theorem,
this knot must be a totally wild knot.
\end{Thm}

From our construction in the previous Theorem, we can also exhibit these knots
as fixed point sets of geometric actions: 

\begin{Cor}
Let $\tau$ be a PL involution of a sphere $S^{n+2}$ whose fixed point set is
a tamely knotted codimension two PL-embedded sphere $S^n$ ($n\neq 3$).  Let $X$ 
be the strict hyperbolization
of the suspension of $S^{n+2}$, $\tau_h$ the induced involution on
$X$, and $\tilde \tau_h$ a lift of the involution to an involution of the 
universal cover 
$\tilde X$.  Finally, let $\tau_\infty$ be the induced involution on the 
boundary at infinity of the universal cover $\tilde X$.  Then the fixed 
subset of $\tau_\infty$ is a knotted, totally wild embedding of $S^{n}$ in 
$S^{n+2}$.  
\end{Cor}

A consequence of the previous corollary is the following:

\begin{Cor}
There exists a smooth closed manifold $M^{n+3}$ ($n\geq 5$), supporting a locally 
CAT(-1) metric,
and having an involution $\sigma \in PL(M^{n+3})$, with the property that
\begin{itemize}
\item {\bf either} $M^{n+3}$ \underline {cannot} support a negatively curved Riemannian 
metric,
\item {\bf or} $\sigma$ is \underline {not} homotopic to an involutive isometry for 
\underline {any} negatively curved Riemannian metric on $M^{n+3}$.
\end{itemize}
\end{Cor}

Unfortunately, in the previous Corollary, we do not know which of the two 
possibilities occur.
Before starting with the proofs, we make a few remarks about our results.

\begin{Rmk}
Note that geometric involutions induce algebraic involutions, by setting
$\Gamma=\pi_1(X)$.  As such, our Corollary 1.1 also applies to geometric 
involutions, and yields that for geometric involutions whose fixed point
sets are $S^n\hookrightarrow S^{n+2}$, tameness implies the embedding
is unknotted.  Conversely, the examples we give in Corollary 1.3 of knotted
fixed point sets for geometric involutions automatically give us knotted
fixed point sets for algebraic involutions.
\end{Rmk}

\begin{Rmk}
Recall that, in high dimension, examples of PL involutions on spheres 
with tamely knotted 
codimension two spheres as fixed point sets exist, by results of Giffen [16].  
In particular,
there are examples satisfying the hypotheses of our Corollary 1.4.
Concrete examples can be obtained by looking at
Brieskorn spheres, see for instance L\"{u} [19].  This should 
be contrasted with the situation in $S^3$, where by the celebrated
solution to the Smith conjecture (see Bass-Morgan [3]), a tamely 
knotted $S^1$ 
cannot be the fixed point set of a smooth $\mZ_p$-action on $S^3$.  
In contrast, our
Corollary 1.1 shows that tamely knotted codimension two spheres cannot 
arise as fixed point sets of algebraic actions.
\end{Rmk}

\begin{Rmk}
Our approach to Theorem 1.2 is based on the general philosophy that properties 
of the links
in a simplicial complex $X$ give corresponding local properties for the boundary
at infinity $Z$ of the universal cover of the hyperbolization of the simplicial
complex.  Various other results that are based on the same philosophy include:
\begin{itemize}
\item if the links of $X$ are spheres, then the boundary $Z$ is a \v{C}ech
cohomology sphere (Bestvina [4]).
\item if the links of $X$ are PL-spheres, then the boundary $Z$ is a topological
sphere (Davis \& Januszkiewicz [10]).
\item bad homology in the links can propagate to give infinitely generated \v{C}ech 
cohomology in the boundary (see Farrell-Lafont [14] for a more precise statement).
\end{itemize}
It would be interesting to construct other types of pathologies via the 
hyperbolization technique.  It would also be interesting to
see if other methods could be used to construct similar pathologies.
\end{Rmk}

\begin{Rmk} This paper was motivated in part by the following more specific
question (which is still open).  Let $\Gamma =\pi_1(M)$ where $M$ is a
closed negatively curved Riemannian manifold, and let $\alpha:\Gamma \rightarrow
\Gamma$ be an automorphism with $\alpha ^2=Id_\Gamma$.  Is $\alpha$ induced by an 
involution of $M$?  That is to 
say, does there exist a self-homeomorphism $f:M\rightarrow M$ with $f^2=Id_M$, and
$f_\sharp =\alpha$?
\end{Rmk}

\section{Non-existence of tame knots.}

We start by giving a proof of Theorem 1.1:

\begin{Prf}[Theorem 1.1]
We have $\Gamma$ a $\delta$-hyperbolic group, and $\Lambda \leq \Gamma$ a 
quasi-convex subgroup.  We are assuming that $\partial ^\infty \Gamma
=S^{n+2}_\infty$, $\partial ^\infty \Lambda=S^n_\infty$, and since $\Lambda$
is a quasi-convex subgroup of $\Gamma$, we get an embedding $S^n_\infty 
\hookrightarrow
S^{n+2}_\infty$.  We start by showing that, if the embedding is tame, then the 
embedding is the trivial knot.  To argue by contradiction, we make the:

\vskip 10pt

\noindent {\bf Assumption:} We assume $S^n_\infty$ is tame in $S^{n+2}_\infty$.

\vskip 10pt

Since the embedding is tame, we let $N$ be a closed tubular
neighborhood of $S^n_\infty$ in $S^{n+2}_\infty$ ($N_\circ$ it's interior), and let 
$W=S^{n+2}_\infty -N$.  Note that $W^{n+2}$ is an open manifold, homeomorphic to
$S^{n+2}_\infty-S^n_\infty$, and the inclusion $W^{n+2}\subset S^{n+2}_\infty-
S^n_\infty$ is a homotopy equivalence.  For a subset $A\subset S^{n+2}_\infty-
S^n_\infty$, we denote by $\w{A}$ the pre-image $p^{-1}(A)$ where $p:\w{S^{n+2}
_\infty-S^n_\infty}\rightarrow S^{n+2}_\infty-S^n_\infty$ is the covering projection.
Since the embedding was assumed to be a non-trivial tame knot, $n\neq 0,3$ and 
Theorem 4.1 in the Appendix tells us that either:
\begin{itemize}
\item $\pi_1(W)\neq \mZ$, or
\item $H_s(\w{W})\neq 0$ for some $0<s$, $s\neq n$.
\end{itemize}

Let us consider the case where $\pi_1(W)\neq \mZ$.  Since the abelianization of
$\pi_1(W)$ is $\mZ$, we know that the group $\pi_1(W)$ is non-abelian.
Note that $\Lambda$ acts on
$(S^{n+2}_\infty,S^n_\infty)$ by homeomorphisms, and the dynamics of the action 
are sink/source (with the sinks and sources lying on $S^n_\infty$).  
Pick $g$ an element of infinite order in $\Lambda$, and 
$n$ sufficiently large so that we have $g^n(W)\subset N_\circ-S^n_\infty$.  Since 
$g$ is a homeomorphism, we have that $g^n(W)\subset S^{n+2}_\infty-S^n_\infty$
is also a homotopy equivalence.  But the inclusion factors through the inclusion
$N_\circ-S^n_\infty\subset S^{n+2}_\infty-S^n_\infty$.  Since $N$ is a tubular
neighborhood of 
$S^n_\infty$, we see that $\pi_1(N_\circ-S^n_\infty)=\pi_1(S^n\times S^1)=\mZ$ 
(if $n>1$) or $\mZ\oplus \mZ$ (if $n=1$), hence is always abelian.  
This yields a commutative diagram:
$$\xymatrix{ \pi_1(g^n(W)) \ar[rr] \ar[dr]_{\cong}
 & & \pi_1(N_\circ-S^n_\infty) \ar[dl]
\\ & \pi_1(S^{n+2}_\infty-S^n_\infty)  &
}
$$
where the arrows are induced by the inclusions.  Since $\pi_1(S^{n+2}_\infty-
S^n_\infty)$ is non-abelian, but the isomorphism factors through an abelian
group, we get a contradiction.

Next we consider the case where $H_s(\w{W})\neq 0$ for some $0<s$, $s\neq n$.
Arguing as before, we take $g$ to be an element of infinite order in $\Lambda$, and
$n$ sufficiently large so that we have $g^n(W)\subset N_\circ-S^n_\infty$.  Since
$g$ is a homeomorphism, it preserves $H_1(S^{n+2}_\infty-S^n_\infty)$, and
lifts to a self-homeomorphism of $\w{S^{n+2}_\infty-S^n_\infty}$ (which we also 
denote by $g$). Note 
that we have $\w{g^n(W)}=g^n(\w{W})$.  Furthermore, we have that $\w{N_\circ-
S^n_\infty}$
is homeomorphic to $S^n\times \mR^2$, forcing $H_s(\w{N_\circ-S^n_\infty})=0$ (since
$s>0, s\neq n$).  We now have a commutative diagram:
$$\xymatrix{ H_s(\w{g^n(W)}) \ar[rr] \ar[dr]_{\cong}
 & & H_s(\w{N_\circ-S^n_\infty})= 0 \ar[dl]
\\ & H_s(\w{S^{n+2}_\infty-S^n_\infty})  &
}$$
where the arrows are induced by the corresponding inclusions.  Since $H_s(\w{W})$
was assumed to be non-zero, we again get a contradiction.

So in all cases, we see that the embedding $S^n_\infty \hookrightarrow 
S^{n+2}_\infty$ cannot be tame.  Now let $p\in S^n_\infty$ be a wild point.  Since 
the orbit of $p$ under the action of $\Lambda$ on the pair $(S^{n+2}_\infty, 
S^n_\infty)$ is dense, we immediately obtain that $S^n_\infty$ cannot contain 
{\it any} tame points, i.e. the embedding $S^n_\infty\hookrightarrow S^{n+2}_\infty$ 
is totally wild.  This completes the proof of Theorem 1.1.
\end{Prf}

\vskip 5pt

\begin{Prf}[Corollary 1.1]
Let $\tau$ be an automorphism of finite order of the $\delta$-hyperbolic group $\Gamma$, 
and 
$\tau_\infty$ the induced homeomorphism on $\partial ^\infty \Gamma =S^{n+2}_\infty$.
Assume that the fixed point set of $\tau_\infty$ is a tamely embedded 
codimension two sphere $S^n_\infty \hookrightarrow S^{n+2}_\infty$.  We need to
show it is unknotted.

First note that Proposition 2.1 in [14] shows that the fixed set of $\tau_\infty$ 
coincides with $\partial ^\infty (\Gamma ^\tau)$.  A result of Neumann [21] ensures
that $\Gamma^\tau$ is a quasi-convex subgroup of $\Gamma$.  Applying Theorem 1.1, 
we see that this forces the fixed point set to be a trivial knot.
\end{Prf}

\vskip 5pt

\begin{Prf}[Corollary 1.2]
This is immediate, since the hypotheses on $X, Y$ ensure that $\Gamma:=\pi_1(X), 
\Lambda:=\pi_1(Y)$ are $\delta$-hyperbolic groups, and that the embedding 
$\Lambda\leq \Gamma$ is quasi-convex.  Since the embedding $\tilde Y\subset \tilde X$ 
is coarsely equivalent to the embedding $\Lambda\leq \Gamma$, we have a natural 
homeomorphism allowing us to identify $(\partial ^\infty \tilde X, \partial ^\infty
\tilde Y)$ with the pair $(\partial ^\infty \Gamma, \partial ^\infty \Lambda)$.  
The corollary now follows from Theorem 1.1
\end{Prf}

\vskip 5pt

\begin{Prf}[Corollary 1.3]
This follows immediately from Corollary 1.2.
\end{Prf}

\section{Existence of wild knots.}

In this section, we start by focusing on proving Theorem 1.2.  A proof of
Corollaries 1.4 and 1.5 will be given at the end of the section.

We start by observing that the links of vertices in both $\Sigma S^n$ and
the corresponding $\Sigma S^{n+2}$ are PL-spheres.  Since hyperbolization
preserves the local structure of the links, we have that the links of 
vertices in both $X^{n+3}$ and $Y^{n+1}$ are PL-homeomorphic to spheres of
dimension $S^{n+2}$ and $S^n$ respectively.  By a result of Davis-Januszkiewicz 
[10], this implies that the boundaries at infinity of 
the universal covers $\tilde X^{n+3}$ and
$\tilde Y^{n+1}$ are homeomorphic to $S^{n+2}$ and $S^n$ respectively.  
We want to show that the inclusion is a wildly embedded codimension two 
sphere.  We use the subscript ``$\infty$'' to distinguish the
spheres at infinity from any of the other spheres that will appear in our 
proof. 

Let $*\in \tilde Y$ be a lift of the
hyperbolization of the (northern) suspension point, and let $S^n_*$, 
$S^{n+2}_*$ denote the links of the 
point $*$ in $\tilde Y$, $\tilde X$ respectively.  
By Alexander duality, both $S^{n+2}_\infty-S^n_\infty$ and 
$S^{n+2}_*-S^n_*$
are homology circles.  Hence there are canonical infinite cyclic covering
spaces:
$$\w{S^{n+2}_\infty-S^n_\infty} \rightarrow {S^{n+2}_\infty-S^n_\infty}$$
$$\w{S^{n+2}_*-S^n_*} \rightarrow {S^{n+2}_*-S^n_*}$$
If $W$ is any subspace of ${S^{n+2}_\infty -S^n_\infty}$ (or of ${S^{n+2}_*-S^n_*}$), we
denote by $\w{W}$ its inverse image under the appropriate covering
projection.  Note that the restriction $\w{W}\rightarrow W$ is also an 
infinite cyclic covering projection. 

Before starting with the argument, we introduce a little more notation: let $K_*$ 
denote an open metric ball of radius $\epsilon$
centered at $*$ in $\tilde X$, with $\epsilon$ small enough that $\partial K_*$ 
be the link $S^{n+2}_*$, and let:
$$\rho: \bar X -K_* \rightarrow \partial K_* = S^{n+2}_*$$
be the geodesic retraction.  

The next Proposition will be needed in order for us to apply Lefschetz duality
to sets of the form $\rho ^{-1}(U)$, and will primarily be used in the 
proof of Proposition 2.

\begin{Prop}
For each open set $U\subset S^{n+2}_*$, the map $\rho: \rho ^{-1}(U)
\rightarrow U$ is a proper homotopy equivalence.
\end{Prop}

\begin{Prf}[Proposition 1]

By a result of Edwards [12], it is sufficient to show that the pre-image $\rho ^{-1}(p)$ 
is contractible for any point $p\in S^{n+2}_*$.  Note that contractibility of the 
set $P:=\rho ^{-1}(p)\cap X$ would be immediate if we knew that it is a totally 
geodesic subspace of $X$.  Indeed, if this
was the case, the fact that the embedding is totally geodesic would imply that 
$\rho^{-1}(p)$ is homeomorphic to the compactification $P \cup \partial
^\infty P$ of $P$ viewed as a CAT(-1) space.  But the latter is automatically 
contractible. 

Unfortunately, $P$ is not quite totally geodesic.  However, it is a 
{\it quasi-convex} subset of $X$, in the sense that there exists a constant
$C$ such that any geodesic segment with endpoints in $P$ lies in a 
$C$-neighborhood of $P$.  To see this, we merely observe that the CAT(-1)
space $X$ is automatically $\delta$-hyperbolic.  Now given any two 
points $x,y\in P$, we consider the geodesic triangle with vertices $p,x,y$.
Note that the geodesics $\overline{xp}$ and $\overline{yp}$ both lie in the set $P$,
so that $\delta$-hyperbolicity tells us that the geodesic $\overline{xy}$ lies
in a $C$-neighborhood of $\overline{xp}\cup \overline{yp}\subset P$ for some 
uniform constant $C$.

Now a quasi-convex subset of a $\delta$-hyperbolic space is still a
$\delta$-hyperbolic space.  Furthermore, the embedding of the quasi-convex subset 
extends to the 
boundary at infinity.  This implies that  $\rho^{-1}(p)$ coincides with
$P \cup \partial ^\infty P$ (where we view $P$ as a $\delta$-hyperbolic space).  
We can now 
homotope $P \cup \partial ^\infty P$ into the set $P$ (since $\partial ^\infty P$
is a Z-set).  Finally, applying the geodesic retraction to $P$, we see that the
set $\bar P$ is indeed contractible, concluding the proof of Proposition 1.
\end{Prf}

To simplify notation, for a set $Z\subset S^{n+2}_*$, we use $Z_\infty$ to 
denote $\rho^{-1}(Z)\cap S^{n+2}_\infty$, i.e. $Z_\infty$ is the subset of
the sphere at infinity that projects to $Z$ under geodesic retraction.  Our
next proposition allows us to relate the topology of an open set $U$ in 
$S^{n+2}_*-S^n_*$ with the topology of the set $U_\infty$.

\begin{Prop}
For each open set $U\subset S^{n+2}_*-S^n_*$, the map $\rho_\infty: U_\infty 
\rightarrow U$ is a proper homotopy equivalence (where $\rho_\infty$ is the 
restriction of the map $\rho$ to the set $U_\infty$).
\end{Prop}

\begin{Prf}[Proposition 2]
We note that $U_\infty$ is the inverse limit of the sequence $\{U_r\}$, where the
$U_r$ are the intersections of $\rho ^{-1}(U)$ with spheres of radius $r$ centered
at the point $*$ (in particular, $U_\epsilon=U$), and the bonding maps are the 
restriction of the geodesic retraction
$\rho_{r,s}:U_r\rightarrow U_s$ (with $\epsilon\leq s\leq r<\infty$).  We denote
by $\rho_{\infty,j}:\invlim \{U_i\}\rightarrow U_j$ the canonical map from the
inverse limit to the individual $U_j$ (and observe that $\rho_{\infty,\epsilon}$ 
coincides with the map $\rho_\infty$).

We say a map is {\it cell-like} provided the pre-image of each point has the
same shape as a point (where shape refers to the functor from the homotopy 
category to the shape category).  Davis-Januszkiewicz have shown that the bonding
maps $\rho_{r,s}$ are cell-like maps (see Section 3 in [10]).  A basic property
of the shape functor is that it commutes with inverse limits (see
Dydak-Segal [11]).  So under
our hypotheses, we see that:
\begin{align*}
Shape(\rho_{\infty,j}^{-1}(x)) & =Shape(\invlim \{\rho_{i,j}^{-1}(x)\})=\invlim 
\{Shape(\rho^{-1}_{i,j}(x))\} \\  
& =\invlim\{Shape(x)\}= Shape(x)
\end{align*}
where $x\in U$ is arbitrary.  This implies that the canonical map $\rho_{\infty,
\epsilon}=\rho_\infty$
from $\invlim \{U_i\}=U_\infty$ to $U$ is cell-like.  Now a result of
Edwards (Section 4 in [12]) asserts that a cell-like proper surjection of ANR's 
is a proper homotopy equivalence, concluding the proof of Proposition 2.
\end{Prf}

From the previous two Propositions, we can immediately obtain:

\begin{Prop}
For each open set $U\subset S^{n+2}_\infty-S^n_\infty$, the inclusion map
$i^+:U_\infty\hookrightarrow \rho^{-1}(U)$ is a proper homotopy equivalence.
\end{Prop}

\begin{Prf}[Proposition 3]
We observe that we have a commutative diagram:
$$\xymatrix{ U_\infty \ar[rr]^{i^+} \ar[dr]_{\rho_\infty}
 & & \rho^{-1}(U) \ar[dl]^{\rho}
\\ & U  &
}
$$
From Propositions 1 and 2, we know that both of the maps $\rho_\infty$ and $\rho$
are proper homotopy equivalences.  This forces $i^+$ to likewise be a 
proper homotopy equivalence.
\end{Prf}

\vskip 5pt

\begin{Rmk}
We point out that an immediate consequence of Proposition 2 is that the map
$\rho_\infty: U_\infty\rightarrow U$ is a {\it near-homeomorphism} (i.e. it
can be approximated arbitrarily closely by homeomorphisms).  This is a consequence
of the fact that $\rho_\infty$ is cell-like, along with the result of Siebenmann
[23] that cell-like maps between manifolds are near-homeomorphisms.  However,
this stronger result will not used in the rest of this paper.

In the paper of Davis-Januszkiewicz [10], the Siebenmann result was used at the 
level of the maps
$\rho_{r,s}$, giving that the boundary at infinity was the inverse limit of 
topological spheres with bonding maps that were near-homeomorphisms.  A result of
Brown [8] implies that the inverse limit of homeomorphic {\it compact} spaces, 
with bonding maps that are near-homeomorphisms, has to be homeomorphic to the 
spaces in question.  We cannot apply Brown's result in our setting, as we are 
working with {\it open} sets $U$.

Finally, we point out that the argument in Proposition 2 can be used to show the
following: if one has an inverse limit of homeomorphic {\it manifolds}, with bonding 
maps which are cell-like, then the inverse limit is homeomorphic to the manifold in 
question.  In the situation where the manifolds are compact, this follows from 
Brown's result.  In the non-compact case, we reach the same conclusion as Brown,
but require the stronger hypothesis on the bonding maps.
\end{Rmk}

\vskip 5pt

Continuing our proof, notice that $(\rho ^{-1}(U); U, U_\infty)$ is
an open cobordism (of dimension $n+3$); i.e.
$$\partial (\rho ^{-1}(U))= U\amalg U_\infty.$$
We will use $i^+:U_\infty\hookrightarrow \rho^{-1}(U)$ and $i^-:U\hookrightarrow
\rho^{-1}(U)$ to denote the respective inclusions.  Note that, by Proposition 1
and Proposition 3, the maps $i^\pm$ are homotopy equivalences.
Now if $U\subset S^{n+2}_*-S^n_*$, then the pullback of the covering space
$\w{U}\rightarrow U$ via $\rho$ gives a canonical infinite cyclic cover
$\w{\rho ^{-1}(U)} \rightarrow \rho ^{-1}(U)$.  In fact, we have a commutative
diagram:

$$\xymatrix{ \w{\rho ^{-1}(U)} \ar[r]^{\w{\rho}} \ar[d]_p & \w{U} \ar[d]
\\ \rho^{-1}(U) \ar[r]_\rho & U
}
$$
where $\w{\rho}$ is the lift of $\rho$.  Note that $\w{\rho}$ is a proper 
homotopy equivalence which is also a retraction.  Duality considerations yield
the following important fact:

\begin{Prop}
$\w{U_\infty} = p^{-1}(U_\infty)$, where $p$ is the covering projection in the 
commutative diagram above.
\end{Prop}

\begin{Prf}[Proposition 4]
Elementary covering space theory reduces the proof of Proposition 4 to showing 
that the 
inclusion map $\sigma ^+: \rho^{-1}(S^{n+2}_*-S^n_*)\cap S^{n+2}_\infty
\subset S^{n+2}_\infty -S^n_\infty$ induces an isomorphism on $H_1(- ; \mZ)$.
For this, we need that $\bar X-K_*$ and $\bar Y -K_*$ are both manifolds, and if 
$n\geq 5$, then by the h-cobordism theorem they are automatically homeomorphic 
to $S^{n+2}\times [0,1]$ and
$S^n\times [0,1]$ respectively.  Set $W^{n+3}=\bar X-(K_*\cup \bar Y)$, and note that
$W$ is an open cobordism with two boundary components: $\partial ^+W=S^{n+2}_\infty
-S^n_\infty$ and $\partial ^-W=S^{n+2}_*-S^n_*$.  We make the:

\vskip 5pt

\noindent {\bf Claim:} Both inclusion maps $\tau ^\pm: \partial ^\pm W\subset W$ induce
isomorphisms on $H_1(- ; \mZ)$.

\vskip 5pt

Assuming the validity of this claim, we now proceed to complete the verification of 
Proposition 4.  Let $\sigma: \rho ^{-1}(\partial ^-W)\subset W$ denote the inclusion 
map, and consider the following commutative diagram of inclusions:

$$\xymatrix{ \rho ^{-1}(\partial ^-W) \ar[r]^\sigma & W 
\\ \partial ^-W \ar[u] \ar[r]_{Id} & \partial ^-W \ar[u]_{\tau^-}
}
$$

Hence Proposition 1 (with $U=\partial ^-W$) and the Claim show that $\sigma_*$ is an 
isomorphism on $H_1(- ;\mZ)$.  Next consider the commutative diagram of inclusions:

$$\xymatrix{ \rho ^{-1}(\partial ^-W) \ar[rrr]^\sigma & & & W 
\\ (\partial ^-W)_\infty \ar[u]^{i^+} \ar[rrr]_{\sigma^+} & & & 
\partial ^+W \ar[u]_{\tau^+}
}
$$

By the Claim, $\tau ^+_*$ is an isomorphism on $H_1$, and we've just shown that 
$\sigma _*$ is an isomorphism.  By Proposition 3, we have that $i^+_*$ is an
isomorphism on $H_1$.  This implies that $\sigma^+_*$ is an isomorphism on $H_1$.

Hence to complete the proof of Proposition 4, it remains to establish the Claim.  
To do this,
it is clearly enough to show that both $H_i(W,\partial ^\pm W)= 0$, for all $i$.
We will do this only for the case $H_i(W, \partial ^-W)$ since the proof of the other
case is completely analogous.  By Lefschetz duality, we have:
$$H_i(W, \partial ^-W)\cong H^{(n+3)-i}_c(W, \partial ^+W).$$
Hence it is equivalent to show that $H^j_c(W, \partial ^+W)= 0$ for all $j$.
And since $H^j_c(W, \partial ^+W)$ is isomorphic to $H^j_c(W-\partial ^+W)$ 
(see below), it suffices
to show that $H^j_c(W-\partial ^+W)= 0$ for all $j$.

Let $X_*$ and $Y_*$ denote $X-\bar K_*$ and $Y-\bar K_*$ respectively (note that 
$X_*$ and $Y_*$
are homeomorphic to $S^{n+2}\times [1,\infty)$ and $S^n\times [1,\infty)$ respectively).
Now consider the exact sequence in cohomology with compact supports for the pair
$(X_*,Y_*)$:
$$H^j_c(X_*)\leftarrow H^j_c(X_*,Y_*)\leftarrow H^{j-1}_c(Y_*).$$
Because of Lefschetz duality
$$H^j_c(X_*) = H^j_c(S^{n+2}\times [1,\infty)) \cong H_{(n+3)-j}
(S^{n+2}\times [1,\infty), S^{n+2}\times 1)= 0.$$
Similarly 
$H^{j-1}_c(Y_*)= 0$, and therefore $H^j_c(X_*,Y_*)= 0$ for all $j$.  Recall
that if $A$ is a closed subspace of $B$ (both locally compact spaces), then
$$H^j_c(B,A)\cong H^j_c(B-A)$$
(note that this was used in the previous paragraph).  Hence 
$$H^j_c(X_*-Y_*)= 0$$
for all $j$.  Since $X_*-Y_*=W-\partial ^+W$, we have completed the proof of Claim,
and hence also of Proposition 4.
\end{Prf}

\vskip 10pt

In other words, the infinite cyclic covering spaces of $U_\infty$ induced from 
${S^{n+2}_\infty-S^n_\infty}$ and ${S^{n+2}_*-S^n_*}$ are consistent. 
An immediate consequence of this is that, for the infinite cyclic coverings,
we have analogues of Propositions 1,2, and 3.  Since we will need the analogue
of Proposition 2, we explicitly state it below.  For any open set 
$U\subset S^{n+2}_*-S^n_*$, associated to the proper
h-cobordism $(\rho^{-1}(U); U, U_\infty)$, we have a canonical infinite cyclic 
covering
by the induced proper h-cobordism $(\w{\rho^{-1}(U)}; \w{U}, \w{U_\infty})$. 

\vskip 5pt

\noindent {\bf Proposition $2^\prime$:}  The map $\w{\rho}:\w{U_\infty}\rightarrow 
\w{U}$ is a proper homotopy equivalence.

\vskip 5pt

Next we note that, since hyperbolization
preserves links, we have by construction 
that $S^n_* \subset S^{n+2}_*$ is a tame knot.  In particular, since $n\neq 3$,
we have (see Theorem 4.1 in Appendix) either:

\begin{itemize}
\item $\pi_1(S^{n+2}_*-S^n_*)\neq \mZ$, or
\item $H_s(\w{S^{n+2}_*-S^n_*})\neq 0$ for some $0<s$ and $s\neq n$.
\end{itemize}

Our approach will consist of
relating the homotopic/homological properties of the embedding
$S^n_\infty\hookrightarrow S^{n+2}_\infty$ to the corresponding properties for
$S^n_*\hookrightarrow S^{n+2}_*$ (which we know to be non-trivially knotted).
We first consider the case
where knottedness is detected by the homology of the canonical infinite cyclic cover, 
by making the additional assumption:

\vskip 10pt

\noindent {\bf Case 1:}  $H_s(\w{S^{n+2}_*-S^n_*})\neq 0$ for some 
$s\neq 0,n$.

\vskip 10pt

Now let $N_1,N_2$ be a pair of closed tubular neighborhoods of $S^n_* 
\subset S^{n+2}_*$ such that $N_1\subset Int(N_2)$.  Note that $N_1$ and
$N_2$ are both homeomorphic to $S^n\times \mD^2$ (since the embedding
$S^n_*\subset S^{n+2}_*$ is tame).  In fact, there is a
homeomorphism taking $S^n_*\subset N_1\subset N_2$ to $S^n\times \{0\}\subset
S^n\times \frac{1}{2}\mD^2 \subset S^n\times \mD^2$ (where $\mD^2$ refers to 
the unit disk in $\mR^2$).  Let $C=S^{n+2}_*-N_1$, and observe that $C$ is
an open codimension zero submanifold of $S^{n+2}_*$ homeomorphic to the 
knot complement $S^{n+2}_*-S^n_*$.  Furthermore, the inclusion $C\subset
S^{n+2}_*-S^n_*$ is a homotopy equivalence.  

\begin{Prop}
The homomorphism:
$$i_*:H_s(\w{C_\infty}) \rightarrow H_s(\w{S^{n+2}_\infty-S^n_\infty})$$
which is induced by the inclusion map
$$i:\w{C_\infty} \rightarrow \w{S^{n+2}_\infty-S^n_\infty}$$
is non-zero.
\end{Prop}

\begin{Prf}[Proposition 5]
Express $S^{n+2}_\infty-S^n_\infty$ as the union of the following two open sets 
$$A=C_\infty$$
$$B=(Int(N_2))_\infty-S^n_\infty .$$
Note that $A\cap B = (Int(N_2)-N_1)_\infty$.  This implies that
$\w{A}\cup \w{B}=\w{S^{n+2}_\infty-S^n_\infty}$ and $\w{A}\cap \w{B}=\w{A\cap B}$, so
applying the Mayer-Vietoris sequence in dimension $s$, we get:
$$H_s(\w{A}\cap \w{B})\rightarrow H_s(\w{A})\oplus H_s(\w{B})\rightarrow
H_s(\w{S^{n+2}_\infty-S^n_\infty}).$$
Hence to verify Proposition 5, it suffices to show $(1)$ $H_s(\w{A})\neq 0$, and
$(2)$ $H_s(\w{A}\cap \w{B})=0$.  

To see $(1)$, apply {\bf Proposition $2^\prime$} to
the set $U=C$.  This yields an isomorphism $H_s(\w{A})\cong H_s(\w{C})$.  But 
since $C\subset S^{n+2}_*-S^n_*$ is a homotopy equivalence, and as we are 
working under the hypothesis of {\bf Case 1}, we have that $H_s(\w{C})\neq 0$,
establishing $(1)$.

To see $(2)$, apply {\bf Proposition $2^\prime$} to the set $U= Int(N_2)-N_1$ 
(note that 
$\w{A}\cap \w{B}= \w{(Int(N_2)-N_1)_\infty}$).  This yields an
isomorphism $H_s(\w{A}\cap \w{B})\cong H_s(\w{Int(N_2)-N_1})$.  But 
$Int(N_2)-N_1$ is homeomorphic to $S^n\times (S^1\times (\frac{1}{2}, 1))$
by the discussion above.  And hence $\w{Int(N_2)-N_1}$ is homeomorphic to 
$S^n\times \mR^2$.  Since $s\neq 0,n$, we have that $H_s(S^n\times \mR^2)=0$,
establishing $(2)$, and completing the proof of Proposition 5.
\end{Prf}

Now Proposition 5 tells us that the embedding $S^n_\infty\hookrightarrow
S^{n+2}_\infty$ is a knotted sphere, concluding the argument for {\bf Case 1}.
We now focus on:

\vskip 10pt

\noindent {\bf Case 2:}  $\pi_1({S^{n+2}_*-S^n_*})\neq \mZ$.

\vskip 10pt

In this case, we claim that $\pi_1(S^{n+2}_\infty-S^n_\infty)$ is likewise
$\neq \mZ$.  In order to see this, consider the sets $N_1$, $N_2$, and $C$ 
defined previously.  We have that $S^n_*\subset N_1\subset N_2$, with each
$N_i$ a closed tubular neighborhood of $S^n_*$ in $S^{n+2}_*$, and $C$ is 
the set $S^{n+2}_*-N_1$.  

\begin{Prop}
The image of the homomorphism:
$$i_\#:\pi_1(C_\infty)\rightarrow \pi_1(S^{n+2}_\infty-S^n_\infty)$$
which is induced by the inclusion map is a non-abelian group.
\end{Prop}

\begin{Prf}[Proposition 6]
 As in the previous Proposition, we decompose $S^{n+2}_\infty-S^n_\infty$ into a 
 pair of open sets:
$$A=C_\infty$$
$$B=(Int(N_2))_\infty -S^n_\infty$$
Note that $D:=A\cap B=(Int(N_2)-N_1)_\infty$, and $A\cup B=S^{n+2}_\infty-S^n_\infty$.
Applying Siefert-Van Kampen, we see that:
$$\pi_1(S^{n+2}_\infty-S^n_\infty)=\pi_1(A)*_{\pi_1(D)}\pi_1(B)$$
In order to show that $i_\#(\pi_1(A))$ is non-abelian, it is sufficient
to have $(1)$ $\pi_1(A)$ is non-abelian, and $(2)$ $i_\#$ is a monomorphism.

Note that, applying Proposition 2 with $U=C$, we get that $\pi_1(A)\cong
\pi_1(C)$.  Since $C\subset S^{n+2}_*-S^n_*$ is a homotopy equivalence, and as 
we are working under the hypothesis of {\bf Case 2}, we have that 
$\pi_1(A)$ is a non-abelian group (since it's abelianization is
$H_1(A)=\mZ$, while $\pi_1(A)\neq \mZ$).  This gives assertion $(1)$.

Likewise, we can apply Proposition 2 with $U=Int(N_2)-N_1$, obtaining an 
isomorphism $\pi_1(D)\cong \pi_1(Int(N_2)-N_1)$.  But 
$Int(N_2)-N_1$ is homeomorphic to $S^n\times (S^1\times (\frac{1}{2}, 1))$
by the discussion above, which implies that $\pi_1(D)=\mZ$.  Hence to get 
assertion $(2)$, it is enough to show that $\pi_1(D)$ injects into 
$\pi_1(S^{n+2}_\infty-S^n_\infty)$.  In order to see this, observe that 
we have the commutative diagram:

$$\xymatrix{ (Int(N_2)-N_1)_\infty \ar[d]_\rho \ar[rr] & & (S^{n+2}_*-S^n_*)
_\infty \ar[d]^\rho\\ 
Int(N_2)-N_1 \ar[rr] & & S^{n+2}_*-S^n_*
}
$$
where the vertical arrows are given by geodesic retraction, and the horizontal
arrows are inclusions.  Applying the $H_1$ functor we obtain:
$$\xymatrix{ H_1((Int(N_2)-N_1)_\infty) \ar[d]^{\rho_*} \ar[rr] & & H_1((S^{n+2}_*
-S^n_*)_\infty) \ar[d]^{\rho_*}\\ 
H_1(Int(N_2)-N_1) \ar[rr] & & H_1(S^{n+2}_*-S^n_*)
}
$$
Proposition 2 tells us that the vertical maps are isomorphisms.  Since the 
inclusion map $Int(N_2)-N_1 \subset S^{n+2}_*-S^n_*$ induces an isomorphism
on $H_1$,
the bottow arrow in the commutative diagram is also an isomorphism.  This yields
that the inclusion $D=(Int(N_2)-N_1)_\infty \subset (S^{n+2}_*-S^n_*)_\infty$
induces an isomorphism on $H_1$.  Finally, we note that we have a commutative
diagram:
$$\xymatrix{ \pi_1(D) \ar[d]^{\cong} \ar[rrr] & & & \pi_1((S^{n+2}_*
-S^n_*)_\infty) \ar[d]^{ab}\\ 
H_1(D) \ar[rrr]^\cong & & & H_1((S^{n+2}_*-S^n_*)_\infty)
}
$$
where the horizontal arrows are induced by the inclusion $D\subset 
(S^{n+2}_*-S^n_*)_\infty$, while the vertical arrows are given by abelianization.
Since $\pi_1(D)=\mZ$, abelianization gives an isomorphism.  By the argument 
above, the inclusion $D\subset (S^{n+2}_*-S^n_*)_\infty$ induces an isomorphism
on $H_1$.  This implies that the inclusion $D\subset (S^{n+2}_*-S^n_*)_\infty$
induces a monomorphism on $\pi_1$, giving us assertion (2), and hence, completing
the proof of the Proposition.
\end{Prf}

Proposition 6 now tells us that $\pi_1(S^{n+2}_\infty-
S^n_\infty)\neq \mZ$, and hence the embedding $S^n_\infty\hookrightarrow 
S^{n+2}_\infty$ is also knotted in {\bf Case 2}.  This completes the proof
of the second case, and hence of Theorem 1.2.

\begin{Prf}[Corollary 1.4]
We start by observing that our hypotheses ensure the existence of a triangulation
of the pair $(S^{n+2},S^n)$ such that the involution $\tau$ is a simplicial map.
Now recall that the strict hyperbolization 
procedure of Charney \& Davis (section 7 in [9]) takes a simplicial 
complex and functorially assigns to it a topological space (in fact,
a union of compact hyperbolic manifolds with corners) that supports a metric 
of strict negative curvature (i.e. a locally CAT(-1) metric).  By functoriality 
we mean that simplicial 
isomorphisms induce isometries of the resulting spaces. Let us apply this 
procedure to the suspension of the sphere $\Sigma S^{n+2}$ (respectively 
$\Sigma S^n$), and call the resulting space $X$ (respectively $Y$).

Observe that the involution $\tau$ induces an involution on $\Sigma S^{n+2}$
with fixed point set $\Sigma S^n$, which by functoriality of the hyperbolization
procedure, yields an isometric involution $\tau_h$ of $X$ with 
fixed point set the totally geodesic subspace $Y$.  Now $\tilde \tau_h$
is a lift of this involution to the universal cover $\tilde X$, and the fixed
set of $\tilde \tau_h$ will be a lift $\tilde Y$ of $Y$.  It is easy to see
that the fixed point set for the induced involution $\tau_\infty$ on the 
boundary at infinity is precisely the embedded $\partial ^\infty \tilde Y
\subset \partial ^\infty \tilde X$.  By our theorem, this is a knotted, totally
wild, codimension two embedding of $S^n$ into $S^{n+2}$.
\end{Prf}

\begin{Prf}[Corollary 1.5]
We note that in the case where we start with a triangulation of a {\it smooth}
manifold, the hyperbolization naturally carries a canonical smooth structure
(see Charney-Davis [9]).  In particular, the space $X$ constructed in the 
previous Corollary is a well-defined smooth closed manifold, equipped with a locally
CAT(-1) metric.  Furthermore, the involution $\tau$ is a PL involution of
$X$.  In the notation of the statement of our Corollary, we are letting 
$M^{n+3}=X$ and $\sigma=\tau$.

We now claim that either $M^{n+3}$ supports no negatively curved Riemannian metric,
or that $\sigma$ is not homotopic to an involutive isometry for any negatively
curved Riemannian metric.  To argue by contradiction, let us assume that there
exists a negatively curved Riemannian metric $g$, an involutive isometry 
$\sigma^\prime:(M^{n+3},g)\rightarrow (M^{n+3},g)$, and a homotopy $\sigma^\prime 
\simeq \sigma$.
Since $\sigma^\prime$ is homotopic to $\sigma$, we have that the induced maps
$\sigma^\prime_\infty$ and $\sigma_\infty$ coincide on 
$\partial ^\infty \tilde M^{n+3}$.  In particular, the fixed point set of $\sigma 
^\prime
_\infty$ is a knotted, totally wild embedding of $S^n_\infty$ in $\partial ^\infty 
\tilde M^{n+3}=S^{n+2}_\infty$.
Since $\sigma ^\prime$ is an isometric involution on $(M^{n+3},g)$, the fixed point
set of $\sigma ^\prime$ is a totally geodesic submanifold $N\subset M^{n+3}$.  Since
$\partial ^\infty N$ must coincide with the fixed point set of the $\sigma ^\prime
_\infty$, the submanifold $N$ must have codimension two.

Now consider the isometric embedding $\tilde N^{n+1}\subset \tilde M^{n+3}$ which 
is left 
invariant under the lift $\tilde \sigma ^\prime$ of $\sigma ^\prime$, and let
$*\in \tilde N^{n+1}$ be an arbitrary point.  Let $K$ denote a sphere of
radius $\epsilon$ centered at $*$.  Since we are in a  
{\it Riemannian} manifold of negative curvature, geodesic retraction gives
an actual {\it homeomorphism} from the knotted pair $(S^{n+2}_\infty,S^n_\infty)$ 
to the pair $(K, K\cap \tilde N^{n+1})$.  Finally, we note that $K$ is 
homeomorphic to $S^{n+2}$,
and that $K\cap \tilde N^{n+1}$ is an unknotted $S^n$ in $K$ (as it is the 
intersection of a small metric sphere with a totally
geodesic submanifold passing through the point $*$).  This gives us a contradiction,
and completes the proof of the Corollary.
\end{Prf}

\section{Appendix: classification of tame knots in dimension $\geq 6$.}

In this appendix, we provide a proof of the classification of tame knots in 
high dimensional spheres.  While this result is probably known to experts, we
were unable to find a reference in print.

\begin{Thm}
Let $S^n\hookrightarrow S^{n+2}$ be a tame knot, and let $W^{n+2}$ be the complement
of a small open tubular neighborhood of $S^n$ in $S^{n+2}$ (note that $W^{n+2}$ is
a compact manifold with boundary).  Suppose that $n\geq 4$ (so the ambient sphere has
dimension $\geq 6$) and that $\pi_1(W)=\mZ$.
Then there exists an integer $s$ satisfying:
\begin{itemize}
\item $0<s$ and $s\neq n$,
\item $H_s(\tilde W)\neq 0$.
\end{itemize}
\end{Thm}

Before starting the proof, we observe that, if $\pi_1(W)=\mZ$, then the covering 
$\w{W}$ of $W$ corresponding to the abelianization of $\pi_1$ coincides with the
universal cover $\tilde W$ of $W$.  

\begin{Prf}
Note first that if $\tilde H_i(\tilde W)=0$ for all $i$, then $\tilde W$ is contractible,
and hence $W$ is a homotopy circle.  Then the fibering theorem of Browder-Levine [7],
as extended to the topological category by the work of Kirby-Seibenmann [18],
can be used to conclude that $S^n\hookrightarrow S^{n+2}$ is unknotted.  This 
contradiction shows that we can assume that $H_i(\tilde W)\neq 0$ for some $0<i\leq
n+2$.  Now suppose that the only such integer $s$ such that $\tilde H_s(\tilde W)\neq 0$ 
is $s=n$.  We will show this assumption leads to a contradiction, completing the 
proof of the theorem.  

Observe that our assumption implies that $\pi_i(\tilde W)=0$ for all $i<n$, and
hence that $W$ is a compact manifold with the $n-1$ homotopy type of $S^1$ (in
fact, a $K(\mZ,1)$ can be constructed by adding $n+1,n+2,\ldots$ cells to $W$).
Consequently, $H^i_c(\tilde W)=H^i(\mZ, \mZ[\mZ])$ for all $i\leq n-1$, and since
$2<n-1$, we conclude that $H^2_c(\tilde W)=0$ since the cohomological dimension of
$\mZ$ is one.  But on the other hand, Lefschetz duality yields that $H_n(\tilde W)
\cong H^2_c(\tilde W, \partial \tilde W)$.  Applying the long exact sequence for
cohomology with compact supports to the pair $(\tilde W, \partial \tilde W)$ yields:

$$0=H^2_c(\tilde W) \leftarrow H_n(\tilde W)\leftarrow H^1_c(\partial \tilde W)
\leftarrow H^1_c(W)$$

Since $\partial W=S^1\times S^n$, the same reasoning given above for $W$ yields
that the inclusion $\partial W\rightarrow W$ is an $(n-1)$-equivalence, and hence
the inclusion map $\partial \tilde W\rightarrow \tilde W$ induces an isomorphism on
the functor $H^1_c(-)$.  This gives the contradiction $H_n(\tilde W)=0$, completing
the proof of the theorem.
\end{Prf}

\begin{Rmk}
We point out that, in the case $n=1,2$, a tame knot $S^n\hookrightarrow S^{n+2}$ 
is non-trivial if and only if $\pi_1(S^{n+2}-S^n)\neq \mZ$.  For $n=1$, this 
is due to Papakyriakopoulos [22], while for $n=2$, this is Thm. 11.7.A in 
Freedman-Quinn [15].  In particular, there are algebraic criterions for 
deciding whether a tame embedding of a codimension two sphere is knotted when
the ambient sphere has dimension $\neq 5$.
\end{Rmk}

\section{Concluding Remarks.}

The examples we constructed in Theorem 1.2 give codimension 2 embeddings of spheres
which are fairly pathological.  A reasonable question is whether these embeddings
can be approximated by ``nice embeddings''.  More precisely, we have the:

\vskip 5pt

\noindent {\bf Question:} Can one approximate the topological embedding $S^n_\infty
\hookrightarrow S^{n+2}_\infty$ by locally flat topological embeddings?

\vskip 5pt

Note that, in codimension $1$, Bing [5], Ancel-Cannon [2], and
Ancel [1] have shown that
topological embeddings can always be approximated by locally flat topological 
embeddings.  On the other hand, in codimension $2$, Matsumoto [20] has 
constructed topological embeddings of closed surfaces in $\mR^4$ which are {\it not} 
approximable by locally 
flat topological embeddings (however, his hypotheses require the surfaces to
have positive genus).  

Observe that in the examples we constructed in Theorem 1.2, it was precisely the
``local knottedness'' of $Y$ in $X$ that allowed us to show that the embedding
of the boundaries at infinity were knotted.  One can ask whether this ``local
knottedness'' is really necessary.

\vskip 5pt

\noindent {\bf Question:} Does there exist a compact,
piecewise hyperbolic, 
locally CAT(-1) manifold $X$, containing a totally geodesic submanifold
$Y$, whose universal covers have boundaries at infinity 
$S^{n+2}_\infty$ and $S^n_\infty$ respectively, and having the properties:
\begin{itemize}
\item at each point $y\in Y$, we have that $lk_Y(y)\hookrightarrow lk_X(y)$ is
unknotted (i.e. $Y$ is a locally flat topological submanifold of $X$).
\item $S^n_\infty\hookrightarrow S^{n+2}_\infty$ is knotted.
\end{itemize}

\vskip 5pt

It would be quite interesting if the answer to this last question were positive.
Finally, it is perhaps interesting to observe a formal analogy of Theorem 1.1 with
a totally unrelated result in algebraic K-theory; i.e. the main result in Farrell [13].
That result asserts, for every ring $R$, that the abelian group $Nil(R)$ is either
zero, or not finitely generated.

\section{Bibliography}

\vskip 5pt

\noindent [1]  Ancel, F.D.  {\it Resolving wild embeddings of codimension-one manifolds
in manifolds of dimension greater than $3$},  Topology Apps.  24 (1986), pp. 13--40.

\vskip 5pt

\noindent [2]  Ancel, F.D. \& Cannon, J.W.  {\it The locally flat approximation of 
cell-like embedding relations},  Ann. Math.  109 (1979), pp. 61--86.

\vskip 5pt

\noindent [3]  Bass, H. \& Morgan, J.W.  {\it The Smith conjecture},  Academic
Press, Orlando, FL.  1984, xv+243 pp.

\vskip 5pt

\noindent [4]  Bestvina, M.  {\it Local homology properties of boundaries of groups},
Michigan Math. J.  43 (1996), pp. 123--139.

\vskip 5pt

\noindent [5]  Bing, R.H.  {\it Approximating surfaces by polyhedral ones},  Ann. Math.
65 (1957), pp. 484--500.

\vskip 5pt

\noindent [6]  Bing, R.H.  {\it The geometric topology of 3-manifolds},  Amer. Math. Soc.,
Providence, RI.  1983, x+238 pp.

\vskip 5pt

\noindent [7]  Browder, M. \& Levine, J.  {\it Fibering manifolds over a circle},  
Comment. Math. Helv.  40 (1966), pp. 153--160.

\vskip 5pt

\noindent [8]  Brown, M.  {\it Some applications of an approximation theorem for
inverse limits},  Proc. Amer. Math. Soc.  11 (1960), pp. 478--483.

\vskip 5pt

\noindent [9]  Charney, R.M. \& Davis, M. W.  {\it Strict hyperbolization},
 Topology  34  (1995),  no. 2, pp. 329--350.

\vskip 5pt

\noindent [10]  Davis, M.W. \&  Januszkiewicz, T. {\it Hyperbolization of polyhedra},
 J. Differential Geom.  34  (1991),  no. 2, pp. 347--388.

\vskip 5pt

\noindent [11]  Dydak, J. \& Segal, J.  {\it Shape theory: an introduction}, 
Springer-Verlag, Berlin.  1978, vi+150 pp.

\vskip 5pt

\noindent [12]  Edwards, R.D.  {The topology of manifolds and cell-like maps}, in
Proceedings of the I.C.M., Helsinki (pp. 111--127).  Acad. Sci. Fennica, Helsinki, 
1980.

\vskip 5pt

\noindent [13]  Farrell, F.T.  {\it The non-finiteness of Nil},  Proc. Amer. Math. Soc.
65  (1977), pp. 215--216.

\vskip 5pt

\noindent [14]  Farrell, F.T. \& Lafont, J.-F.  {\it Finite automorphisms of 
negatively curved Poincar\'e Duality groups},  to appear in Geom. Funct. Anal.

\vskip 5pt

\noindent [15]  Freedman, M.H. \& Quinn, F.  {\it Topology of 4-manifolds},  Princeton
University Press, Princeton, NJ.  1990, viii+259 pp.

\vskip 5pt

\noindent [16]  Giffen, C.H.  {\it The generalized Smith conjecture},  Amer. J.
Math.  88 (1966), pp. 187--198.

\vskip 5pt

\noindent [17]  Kirby, R.C. \& Siebenmann, L.C.  {\it Normal bundles for codimension
$2$ locally flat imbeddings}, in Geometric Topology (Proc. Conf., Park City, Utah,
1974) (pp. 310--324).  Springer, Berlin, 1975.

\vskip 5pt

\noindent [18]  Kirby, R.C. \& Siebenmann, L.C.  {\it Foundational essays on topological
manifolds, smoothings, and triangulations},  Princeton University Press, Princeton,
NJ.  1977, vii+355 pp.

\vskip 5pt

\noindent [19]  L\"{u}, Z.  {\it A note on Brieskorn spheres and the generalized 
Smith conjecture},  Michigan Math. J.  47  (2000), no. 2, pp. 325--333.

\vskip 5pt

\noindent [20]  Matsumoto, Y.  {\it Wild embeddings of piecewise linear manifolds in
codimension two}, in Geometric Topology (pp. 393--428).  Academic Press, New York, NY,  
1979. 

\vskip 5pt

\noindent [21]  Neumann, W.D.  {\it The fixed group of an automorphism of a word hyperbolic
group is rational},  Invent. Math. 110  (1992), no. 1, pp. 147--150. 

\vskip 5pt

\noindent [22]  Papakyriakopoulos, C.D.  {\it On Dehn's Lemma and the asphericity of
knots},  Ann. Math. 66  (1957), pp. 1--26.

\vskip 5pt

\noindent [23]  Siebenmann, L.C.  {\it Approximating cellular maps by homeomorphisms},
Topology  11 (1972), pp. 271--294.

\end{document}